\documentclass[11pt]{article}
\usepackage{graphicx}
\usepackage{epsfig}
\usepackage{amsfonts}
\usepackage{amssymb}
\usepackage{amsthm}
\usepackage{amsmath}		%
\usepackage{pdflscape}
\usepackage{tikz}
\usepackage{float}
\usetikzlibrary{decorations.pathreplacing}
\usetikzlibrary{patterns}
\graphicspath{{../figures/}}
\newtheorem{thm}{Theorem}[section]

\newtheorem{lem}[thm]{Lemma}
\newtheorem{prop}[thm]{Proposition}
\theoremstyle{remark}
\newtheorem{rem}[thm]{Remark}
\theoremstyle{definition}
\newtheorem{defi}[thm]{Definition}

\newcommand{\xx}{\gamma}
\newcommand{\yy}{\delta}
\newcommand{\patt}{\mathbf{Patt}}
\newcommand{\Dec}{\mathbf{Dec}}
\newcommand{\Seq}{\mathbf{Seq}}

\begin{document}

\title{Large Deviations for Permutations Avoiding Monotone Patterns}
\author{Neal Madras \\ Department of Mathematics and Statistics \\
York University \\ 4700 Keele Street  \\ Toronto, Ontario  M3J 1P3 Canada 
\\  {\tt  madras@mathstat.yorku.ca}  \\ and \\ Lerna Pehlivan  \\ Department of Mathematics and Statistics \\
Dalhousie University \\ 6316 Coburg Road \\ Halifax,  Nova Scotia B3H 4R2 Canada  \\
{\tt  lr608779@dal.ca }
 }
\maketitle

\begin{abstract} 
For a given permutation $\tau$, let $P_N^{\tau}$ be the uniform probability 
distribution on the set of $N$-element permutations $\sigma$ that avoid the pattern $\tau$.
For $\tau=\mu_k:=123\cdots k$, we consider $P_N^{\mu_k}(\sigma_I=J)$ where
$I\sim \xx N$ and $J\sim \yy N$ for $\xx,\yy\in (0,1)$.  If $\xx+\yy\neq 1$ then we are in 
the large deviations regime with the probability decaying exponentially, and we 
calculate the limiting value of $P_N^{\mu_k}(\sigma_I=J)^{1/N}$. 
We also observe that for $\tau = \lambda_{k,\ell} := 12\ldots\ell k(k-1)\ldots(\ell+1)$
and $\gamma+\delta<1$,  
the limit of $P_N^{\tau}(\sigma_I=J)^{1/N}$ 
is the same  as for $\tau=\mu_k$. 
\end{abstract}

\section{Introduction and Statement of Results}
  \label{sec-intro}

This paper concerns an aspect of the probabilistic properties of a class of
pattern-avoiding permutations.  As surveyed in the books of B\'{o}na \cite{bona1} and 
Kitaev \cite{kitaev}, pattern avoidance has been of considerable
interest in combinatorial theory, interacting with fields ranging from algebraic
combinatorics to the theory of algorithms.  In the next few paragraphs, we 
give a brief description of the context.

For each positive integer $N$, let ${\cal S}_N$ be the set of all permutations of $1,2,\ldots,N$.
We represent a permutation $\sigma \in {\cal S}_N$  as a string of numbers using the 
one-line notation $\sigma=\sigma_1\ldots\sigma_N$.  We also view $\sigma$ as the function
on $\{1,\ldots,N\}$ that maps $i$ to $\sigma(i)=\sigma_i$. 
The graph of the function $\sigma$ is the set of $N$ points 
$\left\{(i,\sigma_i):i=1,\ldots, N\right\}$ in $\mathbb{Z}^2$.
Given $\tau \in {\cal S}_k$ (with $k\leq N$), we say that 
a permutation $\sigma\in {\cal S}_N$ avoids the pattern $\tau$ 
(or ``$\sigma$ is $\tau$-avoiding'') if there is no $k$-element subsequence 
of $\sigma_1,\ldots,\sigma_N$ having the same relative order as $\tau$. 
(See Section \ref{sec-intro2} for a more formal definition.)
Let ${\cal S}_N(\tau)$ be the 
set of permutations in ${\cal S}_N$ that avoid $\tau$. For example, the permutation
$24153$ is not in ${\cal S}_5(312)$ because it contains the subsequence $413$, which has 
the same relative order as $312$.  In contrast,  the permutation $35421$ has no such 
subsequence, and hence $35421 \in {\cal S}_5(312)$.  

We write $|{\cal A}|$ to denote the number of elements in a set ${\cal A}$. 
Knuth \cite{knuth} proved that  $|{\cal S}_N (\tau)|$ is the same for all 
$\tau \in {\cal S}_3 $  and is equal to the $N$th Catalan number, that is 
$\binom{2N}{N}/(N+1)$ for every $N$. For $\tau \in {\cal S}_k$ with $k \geq 4$, the 
values of $|{\cal S}_N (\tau)|$ depend on the pattern $\tau$ and have been computed 
for only some cases.   For example, Gessel \cite{Gessel} used generating functions to 
show that 
\[
     |{\cal S}_N(1234)|   \;=\;   2\sum_{k=0}^N \binom{2k}{k} \binom{N}{k}^2
        \frac{3k^2+2k+1-N-2kN}{(k+1)^2(k+2)(N-k+1)}  \,.
\]

In 2004 Marcus and Tardos \cite{marcus} proved that 
\[L(\tau) \;:=\; \lim_{N\rightarrow\infty}|{\cal S}_N(\tau)|^{1/N} \, \mbox{exists and is finite for every } 
          \tau %
, \]
thereby confirming the Stanley-Wilf Conjecture that had been open for more than two decades. 
For example, for $k\geq 3$ and $1\leq \ell\leq k-2$, consider the patterns 
\[    \mu_k   \;=\; 123\ldots k   \hspace{8mm}\hbox{and} \hspace{8mm}
    \lambda_{k,\ell}  \;=\; 123\ldots (\ell{-}1)\ell k(k{-}1)\ldots (\ell+1)    \,;
\] 
that is, $\mu_k$ is the increasing pattern of length $k$, and  $\lambda_{k,\ell}$ is 
obtained by reversing the last $k{-}\ell$ elements of $\mu_k$. 
A theorem due to Regev \cite{regev} implies that $L(\mu_k)=(k-1)^2$.
Backelin, West and Xin \cite{backelin} prove that
$\mu_k$ and $\lambda_{k,\ell}$ are Wilf equivalent, i.e.\ that $|{\cal S}_N(\mu_k)|\, = \,
 |{\cal S}_N(\lambda_{k,\ell})|$ for every $N$, which
implies that $L(\lambda_{k,\ell})=(k-1)^2$.
More generally, \cite{backelin} finds a bijection from 
${\cal S}_N(\tau_{1}\ldots \tau_{\ell}(\ell{+}1)\cdots(k{-}1)k)$
to ${\cal S}_N(\tau_{1}\ldots \tau_{\ell}k(k{-}1)\ldots(\ell+1))$ for any $\tau\in {\cal S}_{\ell}$.

Recently, some researchers have taken a probabilistic viewpoint towards investigating 
pattern-avoiding 
permutations, especially for patterns in ${\cal S}_3$. 
They have been concerned  with the configurational
properties of a typical $\tau$-avoiding permutation of length $N$---more precisely, of a 
permutation drawn uniformly at random from the set ${\cal S}_N(\tau)$. 
Accordingly, we shall write $P^{\tau}_N$ to denote the uniform probability distribution
over the set ${\cal S}_N(\tau)$.  
The following result, proven independently by Miner and Pak \cite{pak} and 
by Atapour and Madras \cite{madras}, motivates the present paper.

\begin{thm}
  \label{thm.keq3}
   \cite{madras,pak}
Fix numbers $\xx$ and $\yy$ in $(0,1)$ such that $\xx<1-\yy$.
For each $N$, let $I_N$ and $J_N$ be integers in $[1,N]$ such that
\begin{equation}
   \label{eq.IJN}
     \lim_{N\rightarrow\infty}  \frac{I_N}{N} \;=\; \xx 
   \hspace{5mm}\hbox{and}\hspace{5mm}   
   \lim_{N\rightarrow\infty}  \frac{J_N}{N} \;=\; \yy \,.
\end{equation}
Then
\begin{align}
   \label{eq.limP3.1}
      \lim_{N\rightarrow\infty}P_{N}^{123}(\sigma_{I_N}=J_N)^{1/N}  \; & =\;
     \frac{1}{4} \,G(\xx,1-\yy;1)   \\
 \label{eq.limP3.2}
     &  = \; \lim_{N\rightarrow\infty}P_{N}^{132}(\sigma_{I_N}=J_N)^{1/N},
\end{align}
where we define
\begin{equation}
  \label{eq.Gdef1.3} 
    G(u,v;1) \;:=\; \frac{(u+v)^{(u+v)}(2-u-v)^{2-u-v} }{u^u v^v(1-u)^{(1-u)}(1-v)^{(1-v)} } \,.
\end{equation}
\end{thm}
\noindent
Since  $G(u,v;1)\,<\,4$ whenever $u\neq v$, we see that the probabilities
$P_{N}^{123}(\sigma_{I_N}=J_N)$ and $P_{N}^{132}(\sigma_{I_N}=J_N)$
decay exponentially in $N$ when $\xx<1-\yy$.  Thus, a random $123$-avoiding or
$132$-avoiding permutation is exponentially unlikely to contain any points
$\epsilon N$ below the diagonal $\{(i,N{-}i{+}1):1\leq i \leq N\}$; 
we refer to this as the ``large deviations'' regime.
In the case that $\gamma>1-\delta$, Equation (\ref{eq.limP3.1}) still holds (by symmetry about 
the diagonal), but for $\tau=132$ there is no exponential decay---i.e.\ the limit
in Equation (\ref{eq.limP3.2}) is 1.  
In fact, $P_{N}^{132}(\sigma_{I_N}=J_N)$ is asymptotically proportional to $N^{-3/2}$
(\cite{madras3}, \cite{pak}).   
Madras and Pehlivan \cite{madras3} also examined joint probabilities under $P_N^{132}$, 
proving for example that the probability that graph of $\sigma$ has two 
specified points below the diagonal is of order $N^{- 3}$ (under certain conditions on the points). 
Rizzolo, Hoffman, and Slivken \cite{hoffman} 
proved that for $\tau \in {\cal S}_3$, the shape of a $\tau$-avoiding random  permutation can 
be described by Brownian excursion. Janson \cite{janson} studied the number of 
occurrences of another pattern $\pi$ inside a random $132$-avoiding permutation.

Although patterns of length 3 are amenable to precise probabilistic
results, analogues for longer patterns seem to be much harder.
One reason for this is that for $\tau\in {\cal S}_3$, there are nice bijections from ${\cal S}_N(\tau)$
to the set of Dyck paths of length $2N$, and these
bijections translate various configurational properties of $\tau$-avoiding
permutations into tractable properties of Dyck paths (e.g.\ \cite{hoffman},\cite{madras3}).  
(At a more metaphysical
level:  when the Catalan numbers appear in a problem, nice things happen.)
However, nice bijections are much harder to find for patterns of 
length 4.  Although exact formulas for $|{\cal S}_N(\tau)|$ are known for some
patterns $\tau$ of length 4, their proofs are much more complicated
than for length 3 and do not seem to be useful for investigating
properties of $P_N^{\tau}$.   
In this paper our goal is to extend the large deviation result of 
Theorem \ref{thm.keq3}
to the patterns $\mu_k$ for $k\geq 4$.
In contrast to the proof for $\mu_3$, our derivation of the 
precise large deviations results does not require exact formulas
for finite values of $N$.

\begin{figure}[H] %
		\centering
\includegraphics[width=.48 \textwidth]{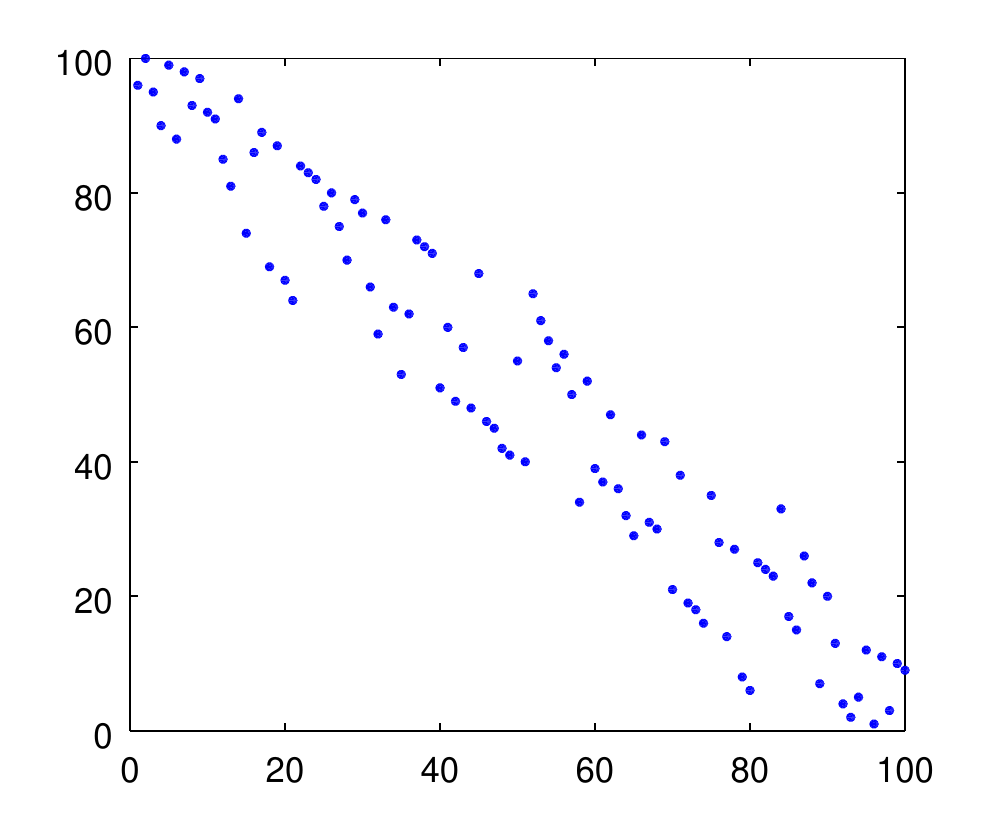}
\includegraphics[width=.48 \textwidth]{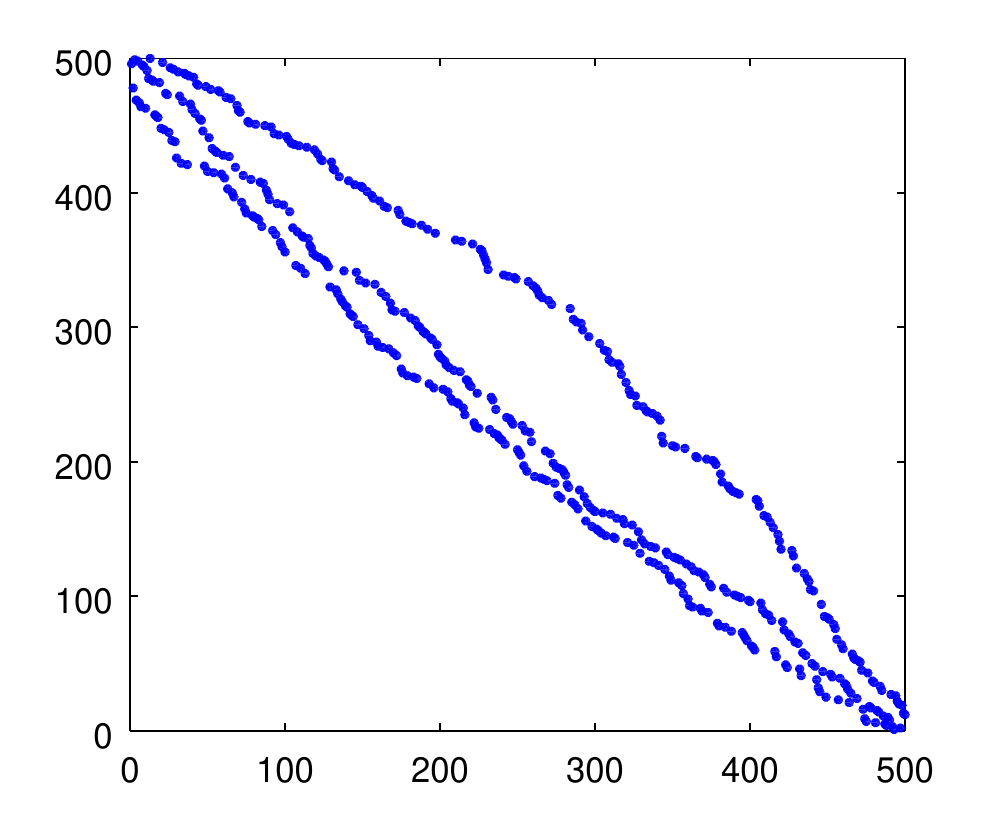}
\caption{Randomly generated 1234-avoiding permutation with N = 100 on the left and N= 500 on the right figure} 
		\label{fig:mc1234}
	\end{figure}

We shall examine the cardinalities of sets of the form
\begin{equation}
   \label{eq.Fdef}
   {\cal F}_{N}(I,J;\tau)   \; := \;   \{\sigma\in {\cal S}_N(\tau):
    \sigma_I=J \} \,.
\end{equation}
Then in terms of the uniform distribution over ${\cal S}_N(\tau)$, we have
\[
    P_N^{\tau}(\sigma_I=J)  \;=\;   \frac{  |{\cal F}(I,J;\tau)|}{|{\cal S}_N(\tau)|}  \,.
\]
Monte Carlo simulations by G\"okhan Y\i ld\i r\i m (as seen in Figure \ref{fig:mc1234}) suggests as $N$ gets larger the number of points well below the $x+y=1$
line decreases. 

We shall typically consider the case $J\ll N-I$ (i.e., points
``below the diagonal''); when $\tau=\mu_k$, the case $J\gg N-I$ follows
from symmetry considerations.  Since we know the asymptotics of 
the denominator $|{\cal S}_N(\tau)|$ for our patterns of interest, and since
our methods are essentially combinatorial, we 
shall henceforth discuss only the numerator,
dealing directly with  $|{\cal F}_N(I,J;\tau)|$ and related combinatorial quantities.

\begin{figure}[H] %
		\centering
		\includegraphics[width=.70 \textwidth,%
		]{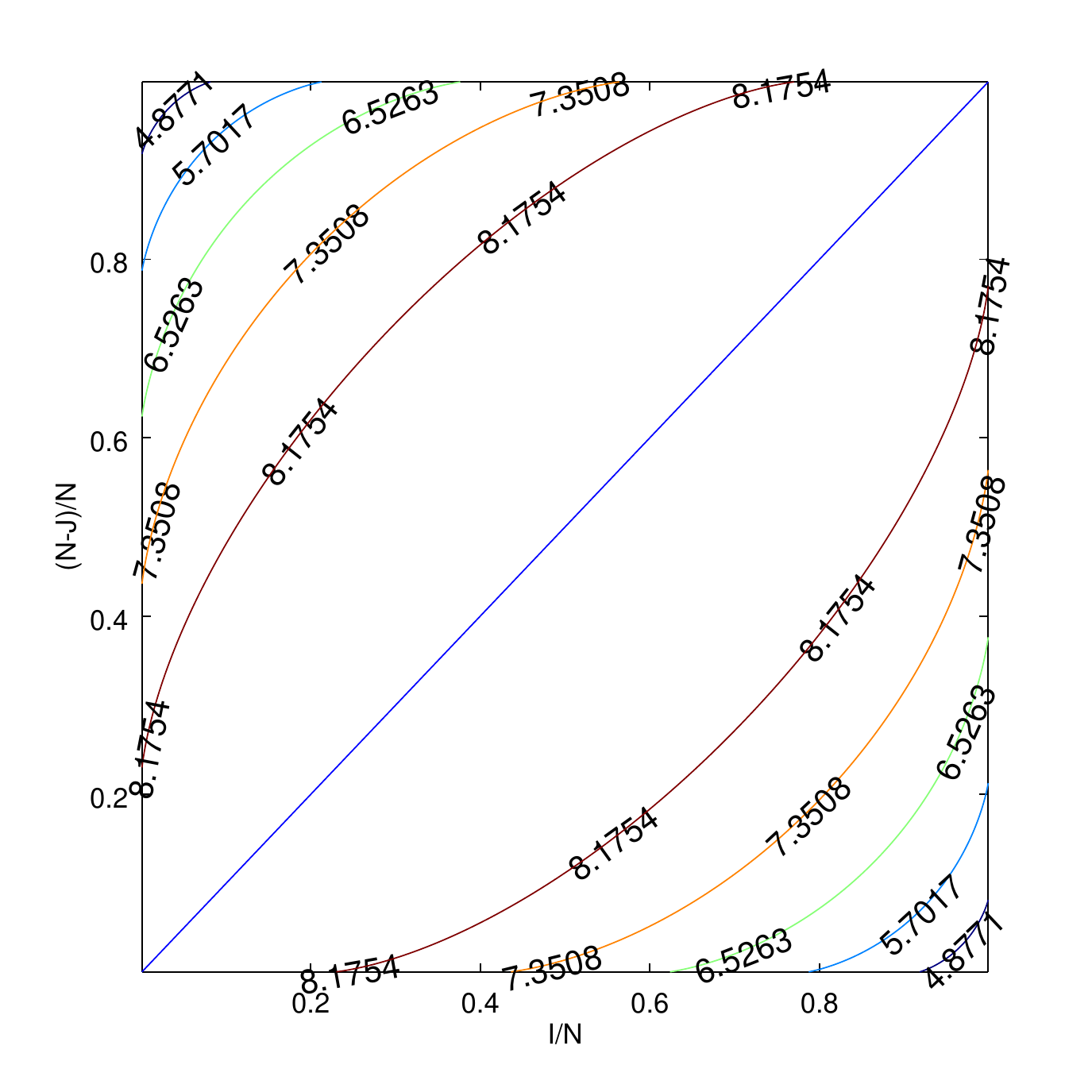}
		\caption{ Level curves in $[0,1]^2$ for $G(u,v;4)$, showing limit values of 
		$|{\cal F}_N(I,J;1234)|^{1/N}$.  To obtain limit values of 
		$P_N^{1234}(\sigma_{I_N}=J_N)^{1/N}$, divide the displayed values by 9.} 
		\label{fig:levelcurves}
	\end{figure}

\begin{thm}
   \label{thm.LD}
Fix $k\geq 4$ and $1\leq \ell\leq k-2$.
Let $\xx$, $\yy$, $I_N$ and $J_N$ be as specified in the statement of Theorem
\ref{thm.keq3}.
Then
\begin{align}
   \label{eq.limF}
      \lim_{N\rightarrow\infty}|{\cal F}_{N}(I_N,J_N;\mu_k)|^{1/N}  \; & =\;
     G(\xx,1-\yy;(k-2)^2)   \\
  \label{eq.limF2}
     & = \; \lim_{N\rightarrow\infty}|{\cal F}_{N}(I_N,J_N;\lambda_{k,\ell})|^{1/N},
\end{align}
where we define
\begin{equation}
  \label{eq.Gdef1} 
    G(u,v;c) \;:=\; 4c\,g(u,v;c)g(v,u;c)g(1-u,1-v;c)g(1-v,1-u;c)
\end{equation}
\begin{equation}
   \label{eq.gdef1}
\hbox{and}   \hspace{10mm}
   g(x,y;c) \;:=\; 
   \left( \frac{2cx+(y-x)-\sqrt{(y-x)^2+4cxy}}{x(c-1)}\right)^{-x} .
\end{equation}
\end{thm}
Figure \ref{fig:levelcurves} gives an example of the level curves of $G(u,v;c)$ for $(u,v) \in [0,1]^2$ and $c=4$. 

\medskip

\bigskip
\begin{rem}
   \label{rem-diag}
When $J_N\approx N-I_N$ (i.e., when we are close to the diagonal), then we 
are in the (limiting) case $\xx=1-\yy$.  
This is not a ``large deviation,''  since $G(u,u;(k-2)^2)=L(\mu_k)$; indeed,
\[   g(x,x;c)  \;=\;   \left(  \frac{2cx -\sqrt{4cx^2}}{x(c-1)} \right)^{-x}    \;=\;
            \left(  \frac{2\sqrt{c}(\sqrt{c}-1)}{c-1} \right)^{-x}   \;=\;
             \left(  \frac{2\sqrt{c}}{\sqrt{c}+1} \right)^{-x} \,,
\]
and it follows that
\[    G(u,u;c) \;=\;  4c\, \left(  \frac{\sqrt{c}+1}{2\sqrt{c}}\right)^2  \;=\;  (\sqrt{c}+1)^2 \,,
\]
which equals $(k-1)^2$ when we substitute $c=(k-2)^2$.
The regime $|N-I_N-J_N|\,=\,o(N)$
is examined by Fineman, Slivken, Rizzolo, and Hoffman (in preparation). 
\end{rem}

\begin{rem}
The numerator and denominator inside the parentheses in Equation 
(\ref{eq.gdef1}) are both 0 when we set $c=1$.  Therefore we 
define $g(x,y;1)$ by taking the limit of $g(x,y;c)$ as $c\rightarrow 1^+$.
We then obtain
\[    g(x,y;1) \;=\; \left(\frac{2x}{x+y}\right)^{-x}
\]
which in turn implies that $G(u,v;1)$ is given by Equation (\ref{eq.Gdef1.3}).
Thus our Theorem \ref{thm.LD} formally recovers Theorem \ref{thm.keq3}.
\end{rem}

\begin{rem}
Assume that $\gamma$, $\delta$, $I_N$ and $J_N$ are as in Theorem \ref{thm.keq3}
except that $\gamma>1-\delta$. Then Equation (\ref{eq.limF}) still holds (by symmetry),
while $\lim_{N\rightarrow\infty}|{\cal F}_{N}(I_N,J_N;\lambda_{k,\ell})|^{1/N}=(k-1)^2$ by 
Proposition 3.1 of \cite{madras} (i.e.,  
$\lim_{N\rightarrow\infty} P_{N}^{\lambda_{k,\ell}}(\sigma_{I_N}=J_N)^{1/N}\,=\,1$).
\end{rem}

\medskip

The term $(k-2)^2$ appears in Equations (\ref{eq.limF}) and (\ref{eq.limF2}) because it is the value of $L(\mu_{k-1})$.
This is highlighted and generalized in Theorem \ref{thm.LDgen} below.

\begin{defi}
  \label{def-SNA}
Let $N$ and $A$ be positive integers, and let $\tau$ be a fixed
permutation.  Define
\[
    {\cal S}_N^{*A}(\tau) \;:=\; \{\sigma \in {\cal S}_N(\tau) \,:\,
      \sigma_i>N-i-A \hbox{ for every $i=1,\ldots,N$ } \} .
\]
Thus, the graph of a permutation in  ${\cal S}_N^{*A}(\tau)$ has no point that is more
than $A$ units below   $\{(i,N+1-i):1\leq i\leq N\}$, the decreasing diagonal  of $[1,N]^2$.
\end{defi}
Then Theorem 1.2 of \cite{madras} implies that for every $\epsilon>0$,
$|{\cal S}_N^{*N\epsilon}(123)|/|{\cal S}_N(123)|$ and 
$|{\cal S}_N^{*N\epsilon}(132)|/|{\cal S}_N(132)|$ converge to 1 
exponentially rapidly as $N\rightarrow\infty$.  

\begin{defi}
   \label{def.oslash}
For $\omega\in {\cal S}_m$,  let $1\oslash\omega$ be the permutation 
$1(\omega_1+1)(\omega_2+1)\ldots(\omega_m+1)$ in ${\cal S}_{m+1}$.
\end{defi}

\noindent
For example, $1\oslash 3124 \,=\,14235$.  Observe that $1\oslash \mu_{k-1}=\mu_k$
and $1\oslash\lambda_{k-1,\ell-1}=\lambda_{k,\ell}$.

Most of the present paper will focus on the proof of the following theorem.

\begin{thm}
   \label{thm.LDgen}
Let $\hat{\tau}$ be  a pattern of length 3 or more, and assume that 
\begin{equation}
    \label{eq.tauhatcond}
       \lim_{N\rightarrow\infty}|{\cal S}_N^{*N\epsilon}(\hat{\tau})|^{1/N}  \;=\; L(\hat{\tau})
        \hspace{5mm}  \hbox{for every $\epsilon>0$}.
\end{equation}
Let $\tau = 1\oslash \hat{\tau}$.  Let $\xx$, $\yy$, $I_N$ and $J_N$ be as 
specified in the statement of Theorem \ref{thm.LD}.
Then
\begin{align}
   \label{eq.limFgen}
      \lim_{N\rightarrow\infty}|{\cal F}_{N}(I_N,J_N;\tau)|^{1/N}  \; & =\;
     G\left(\xx,1-\yy;L(\hat{\tau}) \right).
\end{align}
\end{thm}

\begin{rem}
   \label{rem.LDgen}
(a) Theorem 1.2 of  \cite{madras} implies that Equation (\ref{eq.tauhatcond})
holds for $\mu_3$ and $\lambda_{3,1}$.  
\\
(b) Theorem 1.3(b) of \cite{madras} 
implies that if Equation (\ref{eq.tauhatcond}) holds, then $\hat{\tau}_1$ must equal 1.
The converse of this statement has neither been proved nor disproved; however, simulations
in \cite{madras} and \cite{madras2} suggest that (\ref{eq.tauhatcond}) is false for $\hat{\tau}=1324$. 
\end{rem}

As we shall see in Section \ref{sec-conc}, Theorem \ref{thm.LD} follows from 
Theorem \ref{thm.LDgen} by
induction on $k$, with Remark \ref{rem.LDgen}(a) leading to the base case $k=4$. 
The idea behind the proof of Theorem \ref{thm.LDgen} consists of three  main steps.
An important role is played by the set ${\cal F}_{N}^*(I,J;\tau)$
of permutations in ${\cal F}_{N}(I,J;\tau)$ for which $(I,J)$ is a left-to-right
minimum (i.e., $\sigma_i>J$ for all $i<I$).
The first step is to derive an explicit upper bound to show that $|{\cal F}_{N}^*(I,J;\tau)|^{1/N}$
is less than or equal to $G(\xx,1-\yy; L(\hat{\tau}))$ in the limit.
The second step is to use monotonicity of $G$
to show that we can 
replace ${\cal F}^*$ by ${\cal F}$ in the preceding assertion.
The third step uses the dominant terms from the upper bound of the first step 
to construct a lower bound on 
$|{\cal F}_{N}(I,J;\tau)|^{1/N}$ that is arbitrarily close to the upper
bound.  Section \ref{sec-upperbound} carries out the first two steps,
while Section \ref{sec-lowerbound} performs the third step.
Section \ref{sec-conc} ties the pieces together to complete the proofs of the two
theorems.
Section \ref{sec-intro2} presents some basic definitions and a useful lemma.

We close this section with a physical analogy to help visualize our results about $\mu_k$.
It is easy to verify that an $N$-element permutation $\sigma$ is in ${\cal S}_N(\mu_k)$
if and only if $\sigma$ can be partitioned into $k-1$ decreasing subsequences.
It is not hard to see that these decreasing subsequences are all likely to stay close
to the decreasing diagonal of $[1,N]^2$.  Think of the subsequences as $k{-}1$ elastic
strings, each with one end tied to the point  $(1,N)$ and the other end tied to $(N,1)$, 
and each string tight.   Requiring $\sigma_I$ to equal $J$ is like forcing one of the strings
to pass through the point $(I,J)$.  With this constraint, the rest of the string deforms
into two line segments, one from $(1,N)$ to $(I,J)$ and the other from $(I,J)$
to $(N,1)$.  Tension in the string dictates how the mass of the string is balanced
among the two segments, and  the mass is evenly distributed within each segment.
This physical picture parallels our lower bound construction in Section \ref{sec-lowerbound}.

\subsection{Some Formalities and Preliminaries}
  \label{sec-intro2}

For a string $\omega$ of length $k$ whose entries are all 
distinct numbers, let $\patt(\omega)$ be the
permutation in ${\cal S}_k$ that has the same relative order as $\omega$.
E.g., $\patt(91734) \,=\, 51423$.
More precisely, $\patt(\omega_1\omega_2\cdots\omega_k)$ is the unique
permutation $\pi$ in ${\cal S}_k$ with the property that for all $i,j\in\{1,\ldots,k\}$, 
 $\omega_i<\omega_j$ if and only if $\pi_i<\pi_j$ .

Assume $\tau\in {\cal S}_k$ and $\sigma\in {\cal S}_N$.
We say that $\sigma$ contains the pattern $\tau$ if there exists
$1\leq I_1<I_2<\cdots<I_k\leq N$ such that 
$\patt(\sigma_{I_1}\sigma_{I_2}\cdots\sigma_{I_k})\,=\,\tau$.
We say that $\sigma$ avoids the pattern $\tau$ if $\sigma$ does not contain  $\tau$.
We write ${\cal S}_N(\tau)$ for the set of all permutations in ${\cal S}_N$ that 
avoid $\tau$.

For functions $f$ and $g$, we write $f\sim g$ to mean
$\lim_{N\rightarrow\infty} f(N)/g(N) \,=\,1$.

\begin{defi}
   \label{def.decset}
A finite subset of $\mathbb{Z}^2$ is said to be \textit{decreasing}
if it can be written in the form
$\{(x(m),y(m)):m=1,\ldots,w\}$ with $x(1)<x(2)<\cdots <x(w)$ and
$y(1)>y(2)>\cdots >y(w)$ for some $w\geq 0$.
\end{defi}

We shall also use the following well-known results.

\begin{lem}
   \label{lem-Stir}
(\textit{i}) Let $s$ and $t$ be integers satisfying $0\leq s\leq t$.
Then 
\[   \binom{t}{s}  \;\leq \;   \frac{t^t}{s^s(t-s)^{t-s}}  \,.
\]
(\textit{ii}) Let $\{s_N\}$ and $\{t_N\}$ be sequences of integers with 
$0\leq s_N\leq t_N$ such that $\lim_{N\rightarrow\infty}s_N/N = S$
and $\lim_{N\rightarrow\infty}t_N/N =T$.   Then
\begin{equation}
   \nonumber  %
     \lim_{N\rightarrow\infty}  \binom{t_N}{s_N}^{1/N}   \;=\;   
     \frac{T^T}{S^S(T-S)^{T-S}}  \,.
\end{equation}
In this lemma, we interpret $0^0$ to be 1.
\end{lem}

\noindent
\textbf{Proof}:  Part (\textit{ii}) follows from Stirling's formula,
and part (\textit{i}) is proven for example in Lemma 2.1(b) in \cite{madras}.
\hfill  $\Box$

\section{The Upper Bound}
    \label{sec-upperbound}

We begin with some definitions.
For  a given permutation $\sigma$, define
\begin{equation}
   \label{eq.defL}
    {\cal M}  \; \equiv \;  {\cal M}(\sigma)  \;:=\;   \{(i,\sigma_i):
       \sigma_i<\sigma_t \hbox{ for every $t<i$}\} \,.
\end{equation}
That is, ${\cal M}$ is the set of points of the graph of $\sigma$
corresponding to left-to-right minima.  
Next, let $\sigma\setminus {\cal M}$ be the string consisting of those
$\sigma_t$ such that $(t,\sigma_t)\not\in {\cal M}(\sigma)$.  
Figure \ref{fig:exp} shows an example. 
More generally, if ${\cal A}$ is a subset of $\mathbb{Z}^2$, let 
$\sigma\setminus{\cal A}$ denote the string consisting of those $\sigma_t$
such that $(t,\sigma_t)\not\in {\cal A}$.

\begin{center}
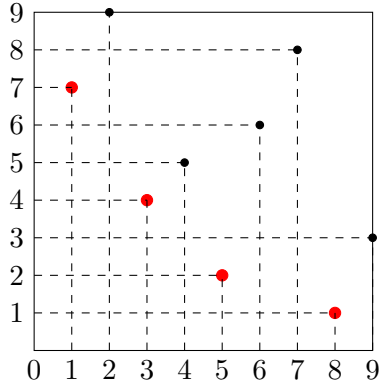
\begin{figure}[htp]
\centering
  \begin{tikzpicture}[scale=0.5]
    \draw (0,0) rectangle (9,9);
    \draw [fill,red] (1,7) circle (0.15);
    \draw [fill]     (2,9) circle (0.1);
    \draw [fill,red] (3,4) circle (0.15);
    \draw [fill]     (4,5) circle (0.1);
    \draw [fill,red] (5,2) circle (0.15);
    \draw [fill]     (6,6) circle (0.1);
		\draw [fill]     (7,8) circle (0.1);
		\draw [fill,red] (8,1) circle (0.15);
		\draw [fill]     (9,3) circle (0.1);

    \draw [dashed] (1,0) -- (1,7);
    \draw [dashed] (2,0) -- (2,9);
    \draw [dashed] (3,0) -- (3,4);
    \draw [dashed] (4,0) -- (4,5);
    \draw [dashed] (5,0) -- (5,2);
		\draw [dashed] (6,0) -- (6,6);
		\draw [dashed] (7,0) -- (7,8);
		\draw [dashed] (8,0) -- (8,1);
		\draw [dashed] (9,0) -- (9,3);

		\draw [dashed] (0,8) -- (7,8);
		\draw [dashed] (0,7) -- (1,7);
		\draw [dashed] (0,6) -- (6,6);
    \draw [dashed] (0,4) -- (3,4);
    \draw [dashed] (0,3) -- (9,3);
    \draw [dashed] (0,5) -- (4,5);
    \draw [dashed] (0,2) -- (5,2);
    \draw [dashed] (0,1) -- (8,1);

    \node [below] at (0,0) {$0$};
    \node [below] at (1,0) {$1$};
    \node [below] at (2,0) {$2$};
    \node [below] at (3,0) {$3$};
    \node [below] at (4,0) {$4$};
    \node [below] at (5,0) {$5$};
    \node [below] at (6,0) {$6$};
    \node [below] at (7,0) {$7$};
		\node [below] at (8,0) {$8$};
		\node [below] at (9,0) {$9$};
		\node [left] at (0,1) {$1$};
    \node [left] at (0,2) {$2$};
    \node [left] at (0,3) {$3$};
    \node [left] at (0,4) {$4$};
    \node [left] at (0,5) {$5$};
    \node [left] at (0,6) {$6$};
		\node [left] at (0,7) {$7$};
		\node [left] at (0,8) {$8$};
		\node [left] at (0,9) {$9$};
  \end{tikzpicture}
  \caption{The graph of $\sigma=794526813 \in S_9$. Here, ${\cal M}=\{(1,7), (3,4), (5,2),(8,1)\}$
and $\sigma\setminus {\cal M} \,=\, 95683$.}
	\label{fig:exp}
  \end{figure}
  \end{center}

The following observations are useful.  We omit the straightforward proof.

\begin{lem}
   \label{lem-minusL}
(\textit{i}) 
A permutation $\sigma$ is uniquely determined by the 
set ${\cal M}$ and the permutation $\patt(\sigma\setminus {\cal M})$.
\\
(\textit{ii}) Let $\hat{\tau}$ be a pattern with $\hat{\tau}_1=1$.
The permutation $\sigma$ avoids $1\oslash \hat{\tau}$ if and only if 
$\patt(\sigma\setminus {\cal M})$ avoids $\hat{\tau}$.
\end{lem}

\medskip
Recall from Section \ref{sec-intro} that
\[   {\cal F}_{N}^*(I,J;\tau) \;=\; \{\sigma\in  {\cal F}_{N}(I,J;\tau) \,:\; \sigma_i>J 
    \hbox{ for all $i<I$}  \,\}.
\]
We shall now perform the first step in the proof of our
main theorem.

\begin{prop}
   \label{prop.Fstar}
Let $\hat{\tau}$ be a pattern of length 3 or more such that $\hat{\tau}_1=1$, and 
let $\tau =1\oslash \hat{\tau}$.
Let $\xx$, $\yy$, $I_N$ and $J_N$ be as 
specified in the statement of Theorem \ref{thm.LD}.
Then
\begin{equation}
   \label{eq.limFstar}
   \limsup_{N\rightarrow\infty}|{\cal F}^*_{N}(I_N,J_N;\tau)|^{1/N}  \;\leq \;
     G(\xx,1-\yy;L(\hat{\tau})).
\end{equation}
\end{prop}

\noindent
\textbf{Proof}:
For $I\in [1,N]$ and $\sigma\in {\cal S}_N$,  we define
\[
   {\cal M}^{<I} \,=\; \{(i,\sigma_i)\in {\cal M}(\sigma): i<I\}
     \hspace{5mm}\hbox{and}\hspace{5mm}
   {\cal M}^{>I} \,=\; \{(i,\sigma_i)\in {\cal M}(\sigma): i>I\} \,.
\]

Fix $I$ and $J$ in $[1,N]$ with $J<N-I$.
Suppose we know that $\sigma\in {\cal F}^*_{N}(I,J;\tau)$, 
$l=|{\cal M}^{<I}|$ and $m=|{\cal M}^{>I}|$.
Then ${\cal M}^{<I}$ is a set of $l$ integral points in
$[1,I)\times(J,N]$, and this set must be decreasing (recall Definition \ref{def.decset}).
Therefore there are at most
$\binom{I-1}{l}\binom{N-J}{l}$ possible realizations of ${\cal M}^{<I}$.
Similarly, there are at most 
$\binom{N-I}{m}\binom{J-1}{m}$ possibilities for ${\cal M}^{>I}$.
Recalling Lemma \ref{lem-minusL}, we obtain the following bound:
\begin{align}
   \nonumber
   |{\cal F}^*_{N} & (I,J;\tau)|  
   \\  \nonumber
   & \leq \;   
     \sum_{l=0}^{I-1} \sum_{m=0}^{J-1} \binom{I{-}1}{l}\binom{N{-}J}{l}
       \binom{N{-}I}{m}\binom{J{-}1}{m} |{\cal S}_{N-l-m-1}(\hat{\tau})|
   \\
   & \leq \; 
   H(I{-}1,N{-}J;L(\hat{\tau}))\,H(J{-}1,N{-}I;L(\hat{\tau})) \,L(\hat{\tau})^{N-1}
   \label{eq.FHH}
\end{align}
where we define
\begin{equation}
   \label{eq.Hdef}
   H(a,b;c)  \; :=\; \sum_{n=0}^a \binom{a}{n}\binom{b}{n}c^{-n} \,.
\end{equation}
In the last step, the bound $|{\cal S}_{N-l-m-1}(\hat{\tau})|
\,\leq \, L(\hat{\tau})^{N-l-m-1}$ is proven in Theorem 1 in \cite{arratia}.

We now wish to bound $H(a,b;c)$ for $a\leq b$ and $c>1$.  
By Lemma \ref{lem-Stir}(\textit{i}), we have
\begin{equation}
   \label{eq.Hfmax}
    H(a,b;c)   \;\leq \;  (a+1) \,\sup\{ f(y;a,b,c): 0\leq y\leq a\}
\end{equation}
where 
\begin{equation}
   \label{eq.fdef}
     f(y;a,b,c)  \;=\; \left(\frac{y}{a}\right)^{-y}
    \left( 1-\frac{y}{a}\right)^{y-a} 
     \left(\frac{y}{b}\right)^{-y}
    \left( 1-\frac{y}{b}\right)^{y-b} c^{-y} \,.
\end{equation}
We now pause to state and prove a lemma, which will also be useful later.

\begin{lem}
   \label{lem-fmax}
Fix real numbers $a,b>0$ and $c>1$.  Define the function $f$ as in 
Equation (\ref{eq.fdef}) for real $y$ in the interval $[0,a\wedge b]$
(where $a\wedge b$ is the minimum of $a$ and $b$).  We interpret $0^0=1$,
which makes $f$ continuous on this interval.
Then there is a unique point $y^*\equiv y^*[a,b,c]$ that maximizes
$f$ in this interval, and $0<y^*<a\wedge b$.  Furthermore,
\begin{equation}
   \label{eq.yformula}
   y^*[a,b,c] \;=\; \frac{\sqrt{(a-b)^2+4cab}-(a+b)}{2(c-1)} 
\end{equation}
and the maximum value of $f$ is 
\begin{equation}
   \label{eq.fmax}
   f(y^*[a,b,c];a,b,c)  \;=\;  2^{a+b}g(a,b;c)\,g(b,a;c) \,,
\end{equation}
where $g$ was defined in Equation (\ref{eq.gdef1}).
\end{lem}

\noindent
\textbf{Proof of Lemma \ref{lem-fmax}:}
By calculus, it is easy to see that $\log f$ is a strictly concave
function of $y$ on $[0,a\wedge b]$, and is maximized at the (unique) point
$y^*\equiv y^*[a,b,c]$ in $(0,a\wedge b)$ that satisfies the equation
\begin{equation}
   \label{eq.yroot}
         (a-y^*)(b-y^*) \;=\; c(y^*)^2 \,.
\end{equation}
Thus Equation (\ref{eq.fdef}) becomes
\begin{align}    
   \nonumber
 f(y^*;a,b,c) \; & =\;
   \frac{a^a\,b^b}{(y^*)^{2y^*}(a-y^*)^{a-y^*}(b-y^*)^{b-y^*}c^{y^*} }
    \\
  \label{eq.1-y}
   & = \;
 \left(1-\frac{y^*}{a}\right)^{-a} \left(1-\frac{y^*}{b}\right)^{-b} 
   \hspace{5mm}\textrm{(using (\ref{eq.yroot}))}.
\end{align}
Solving the quadratic equation (\ref{eq.yroot})  for the positive root
gives
\begin{equation}
  y^*[a,b,c] \;=\; \frac{\sqrt{(a+b)^2+4(c-1)ab}-(a+b)}{2(c-1)} \,,
\end{equation}
which leads to Equation (\ref{eq.yformula}).  Finally, inserting 
(\ref{eq.yformula}) into (\ref{eq.1-y}) gives (\ref{eq.fmax}).
\hfill$\Box$

\medskip

We now return to the proof of Proposition  \ref{prop.Fstar}.
By Equation (\ref{eq.Hfmax}) and  Lemma \ref{lem-fmax}, we have
\begin{equation}
   \label{eq.Hggineq1}
    H(I-1,N-J;c)  \;\leq\;  I\, 2^{N+I-J-1}g(I-1,N-J;c)\,g(N-J,I-1;c) \,.
\end{equation}
By Equation (\ref{eq.IJN}) and the explicit form of $g$, we can 
take the limit in Equation (\ref{eq.Hggineq1}) to get 
\[
   \limsup_{N\rightarrow\infty}H(I_N{-}1,N{-}J_N;c)^{1/N}  \; \leq \;
     2^{1+\xx-\yy}\,g(\xx,1-\yy;c) \,g(1-\yy,\xx;c)  \,.
\]
Similarly, we have
\[
   \limsup_{N\rightarrow\infty}H(J_N{-}1,N{-}I_N;c)^{1/N}  \; \leq \;
      2^{1+\yy-\xx}\,g(\yy,1-\xx;c)
     \,g(1-\xx,\yy;c)  \,.
\]
Proposition \ref{prop.Fstar} now follows directly from the above (with $c=L(\hat{\tau})$) 
and Equation (\ref{eq.FHH}).
\hfill$\Box$

\bigskip

Our next task is to replace ${\cal F}^*_N$ by ${\cal F}_N$ in the statement of 
Proposition \ref{prop.Fstar}.  We shall do this by proving a monotonicity property
of $G$ (Lemma \ref{lem.Gmono}) and then using a compactness argument.

\begin{prop}
   \label{prop.Fnostar}
Under the hypotheses of Proposition \ref{prop.Fstar}, we have
\begin{equation}
   \label{eq.limFnostar}
   \limsup_{N\rightarrow\infty}|{\cal F}_{N}(I_N,J_N;\tau)|^{1/N}  \;\leq \;
     G(\xx,1-\yy;L(\hat{\tau})).
\end{equation}
\end{prop}

We begin by showing that $G$ decreases as we move away from the diagonal.
We emphasize that in this lemma, ``increasing'' and ``decreasing'' are used in their
strict sense.

\begin{lem}
   \label{lem.Gmono}
Fix $c>1$.  The function $G(u,v;c)$ defined in Equation (\ref{eq.Gdef1})
is increasing in $u$ and decreasing in $v$ for $0<u<v<1$.
By symmetry, it is also increasing in $v$ and decreasing in $u$ for
$0<v<u<1$.  In particular, $G$ is maximized when $u=v$, where we 
have 
\begin{equation}
   \label{eq.Guu}
    G(u,u;c) \;=\;  ( \sqrt{c}+1)^2
     \hspace{5mm}\hbox{for every $u\in (0,1)$}.
\end{equation}
\end{lem}

\noindent
\textbf{Proof:}
Recall that Equation (\ref{eq.Guu}) was proved in Remark \ref{rem-diag}.

Since $c$ is fixed, we shall suppress it in the following notation.
Let  $r(u,v) = \sqrt{(v-u)^2+4cuv}$ and
$h(u,v) = [2cu+(v-u)-r(u,v)]/u$.  Then
\begin{align*}
 G(u,v;c)=4c(c-1)^2h(u,v)^{-u}h(v,u)^{-v}h(1{-}u,1{-}v)^{1-u}h(1{-}v,1{-}u)^{1-v}
\end{align*}
and hence
\begin{align}
   \nonumber
\ln G(u,v;c) = & \ln(4c(c-1)^2) - u \ln(h(u,v)) - v \ln(h(v,u)) \\
    \label{eq.lnG}
             & - (1 - u) \ln(h(1{-}u,1{-}v)) - (1 - v) \ln(h(1{-}v,1{-}u)) \,.
\end{align}
By routine calculus and some algebraic manipulation, we obtain
\begin{equation}
   \label{eq.partialh}
  \frac{\partial}{\partial u}\ln(h(u,v))  \;=\;  \frac{v}{u\,r(u,v)}
   \hspace{5mm}\hbox{and} \hspace{5mm}
  \frac{\partial}{\partial u}\ln(h(v,u))  \;=\;  -\,\frac{1}{r(u,v)}  \,.
\end{equation}
Using this and Equation (\ref{eq.lnG}), we can show that
\begin{equation}
   \label{eq.partG}
   \frac{\partial}{\partial u}\ln G(u,v;c) \;= \; - \ln (h(u,v))
   + \ln (h(1{-}u,1{-}v)) \,.
\end{equation}
From this and Equation (\ref{eq.partialh}), we also obtain
\[
\frac{\partial^2}{\partial u^2}\ln G(u,v;c) \; =\; 
 - \,\frac{v}{u \,r(u,v)} \;-\; \frac{(1-v)}{(1-u)\, r(1{-}u,1{-}v)}
   \;<\; 0
\]
for every $u$ and $v$ in $(0,1)$.
Therefore $G(u,v;c)$ is strictly concave in $u$ for fixed $v$
(and, by symmetry, it is strictly concave in $v$ for fixed $u$).

Since $h(u,u)=2c-2\sqrt{c}$ for every $u$, it follows that the 
partial derivative in Equation (\ref{eq.partG}) is zero whenever 
$u=v$.  By symmetry, the same is true for the partial derivative
with respect to $v$.
Combining this with the concavity result of the previous paragraph
completes the proof of the lemma.
\hfill $\Box$

\bigskip
\noindent
\textbf{Proof of Proposition \ref{prop.Fnostar}}:   
It is easy to see that 
$${\cal F}_{N}  (I_N,J_N;\tau) \;  \subseteq  \; 
\bigcup_{1 \leq u \leq I_N, 1\leq t \leq J_N} {\cal F}^*_{N}  (u,t;\tau).$$ 
Let $u(N)$ and $t(N)$ be the values of $u$ and $t$ that maximize
$| {\cal F}^*_{N}(u,t;\tau) |$  over $u$ in $[1,I_N]$ and $t$ in $[1,J_N]$.
Then we have
\begin{equation}
   \label{eq.FN2F}   
     |{\cal F}_{N}  (I_N,J_N;\tau) |  \;\leq \; N^2\, |{\cal F}^*_{N}(u(N),t(N);\tau) | \,.
\end{equation}
Let $LS = \limsup_{N\rightarrow\infty}   |{\cal F}^*_{N}(u(N),t(N);\tau)|^{1/N}$. 
There exists a subsequence ${N'}$ such that $|{\cal F}^*_{N'}(u(N'),t(N');\tau)|^{1/N'}$ 
converges to $LS$.
By compactness of $[0,1]^2$, this subsequence
has a sub-subsequence ${N''}$ for which $(u(N'')/N'', t(N'')/N'')$
converges to a point $(\tilde{u},\tilde{t})$ in $[0,\gamma] \times [0,\delta]$.
Thus Proposition \ref{prop.Fstar} tells us that 
$LS \,\leq \,G(\tilde{u},1-\tilde{t}; L(\hat{\tau}))$.
The monotonicity of $G$ (in Lemma \ref{lem.Gmono}) implies that  
$G(\tilde{u},1-\tilde{t}; L(\hat{\tau})) \leq G(\gamma, 1-\delta; L(\hat{\tau}))$ .
Therefore $LS \leq G(\gamma, 1-\delta; L(\hat{\tau})) $.
Hence, using Equation (\ref{eq.FN2F}), we obtain
 Equation (\ref{eq.limFnostar}).
 \hfill $\Box$

\section{The Lower Bound}
   \label{sec-lowerbound}

To get the lower bound on $|{\cal F}_{N}(I,J;\tau)|$, we shall perform
an explicit construction of some permutations in ${\cal F}^*_{N}(I,J;\tau)$
(this is done in the proof of Proposition \ref{prop-constr} below).
The construction is motivated by examining the dominant terms in our
proof of the upper bound, and showing that they are approximately achieved.

The main result of this section is the following.

\begin{prop}
   \label{prop.Fstar2}
Under the hypotheses of Theorem \ref{thm.LDgen}, we have
\begin{equation}
   \label{eq.limFstar2}
   \liminf_{N\rightarrow\infty}|{\cal F}^*_{N}(I_N,J_N;\tau)|^{1/N}  \;\geq \;
     G(\xx,1-\yy;L(\hat{\tau})).
\end{equation}
\end{prop}

\medskip

 The proof of Proposition \ref{prop.Fstar2}
relies on Proposition \ref{prop-constr} and Lemma \ref{lem-Dec}.
We shall first state these two auxiliary results, then prove 
Proposition \ref{prop.Fstar2}, and conclude the section by proving the two 
auxiliary results.

The construction of Proposition  \ref{prop-constr} uses a positive parameter $A$, 
which will afterwards be of 
the order $N\epsilon$ for fixed small $\epsilon$.  We start with a definition.

\begin{defi}
  \label{def-dec}
Let $w$, $M_1$, and $M_2$ be positive integers, with $w\leq M_1\wedge M_2$.
\\
$\bullet$
Let $\Dec(w;M_1,M_2)$ be the collection of all $w$-element decreasing 
subsets of $\{1,\ldots,M_1\}\times \{1,\ldots,M_2\}$.
(Recall Definition \ref{def.decset}.)
\\
$\bullet$
For given $A>0$, let $\Dec^{*A}(w;M_1,M_2)$ be the collections of 
all $w$-element sets ${\cal B}\in \Dec(w;M_1,M_2)$ such that 
\begin{equation}
   \label{def-DecA}
    y\;<\;M_2\,-\,x\,\frac{M_2}{M_1}\,+\,A    \hspace{5mm}
     \textrm{for all $(x,y)\in{\cal B}$}.
\end{equation}
The collections $\Dec(0;M_1,M_2)$ and $\Dec^{*A}(0;M_1,M_2)$ each 
contain one member:  the empty set.
\end{defi}

Observe that the line $y=M_2-xM_2/M_1$ is the decreasing diagonal
of the rectangle $[0,M_1]\times[0,M_2]$.  Thus, $\Dec^{*A}(w;M_1,M_2)$
is the collection of sets in $\Dec(w;M_1,M_2)$ that rise less
than $A$ above the diagonal.

\begin{prop}
   \label{prop-constr}
Let $\hat{\tau}$ be a pattern of length 3 or more such that $\hat{\tau}_1=1$, and 
let $\tau =1\oslash \hat{\tau}$.
Let $N$, $I$, $J$, and $A$ be positive integers with $J<N-I-2A$.
Let $w_1$ and $w_2$ be integers with 
\begin{equation}
   \label{eq.win}
      0 \,\leq\, w_1 \,\leq \, I-1  
    \hspace{4mm}\hbox{and} \hspace{4mm}
    0\,\leq \,  w_2\,\leq \, J-1\,.
\end{equation}
Then (recall Definitions \ref{def-SNA} and \ref{def-dec})
\begin{align}
   \nonumber
    |{\cal F}^*_{N}(I,J;\tau)|  \;  \geq \;  &
        |\Dec^{*A}(w_1;I-1,N-2A-J)|  \\
   \nonumber
        &   \hspace{5mm} \times\,|\Dec^{*A}(w_2;N-2A-I,J-1)|  \\
   \label{eq.FDDS}
       & \hspace{15mm}\times \, |{\cal S}^{*A}_{N-w_1-w_2-1}(\hat{\tau})|  \,.
\end{align}
\end{prop}

\begin{lem}
   \label{lem-Dec}
Consider sequences of positive integers $w(N)$, $M_1(N)$, $M_2(N)$, and 
$A_N$ such that 
\[
    \lim_{N\rightarrow\infty}\frac{w(N)}{N}\,=\,\theta, \hspace{3mm}
    \lim_{N\rightarrow\infty}\frac{M_1(N)}{N}\,=\,\alpha, \hspace{3mm}
    \lim_{N\rightarrow\infty}\frac{M_2(N)}{N}\,=\,\beta, \hspace{3mm}
    \lim_{N\rightarrow\infty}\frac{A_N}{N}\,=\,\epsilon,
\]
with $0<\theta <\alpha\wedge\beta$ and $\epsilon>0$.   Then
\begin{equation}
   \label{eq.lemDec1}
    \lim_{N\rightarrow\infty}  \frac{|\Dec^{*A_N}(w(N);M_1(N),M_2(N))|}{
       |\Dec(w(N);M_1(N),M_2(N))|}   \;=\;  1    
\end{equation}
and (for $f$ defined by Equation (\ref{eq.fdef}))
\begin{align}
   \label{eq.lemDec2}
     f(\theta;\alpha,\beta,c)\,c^{\theta}   \;& = \;
    \lim_{N\rightarrow\infty}  |\Dec(w(N);M_1(N),M_2(N))|^{1/N}   \\      
   \label{eq.lemDec2A}
    & =\;
    \lim_{N\rightarrow\infty}  |\Dec^{*A_N}(w(N);M_1(N),M_2(N))|^{1/N} 
\end{align}
for any $c$.  (Notice that $f(\theta;\alpha,\beta,c)c^{\theta}$ is independent of $c$ by definition.)
\end{lem}

\medskip

\noindent
\textbf{Proof of Proposition \ref{prop.Fstar2}:}   Let $c=L(\hat{\tau})$.
Choose $\epsilon>0$ such that $\xx<1-\yy-2\epsilon$.  Let $\{A_N\}$ be a 
sequence of positive integers such that 
$\lim_{N\rightarrow\infty}A_N/N\,=\,\epsilon$.
Therefore $J_N<N-I_N-2A_N$ holds for all sufficiently large $N$.

Let $\{w_1(N)\}$ and $\{w_2(N)\}$ be sequences of positive integers such that
\begin{align*}
    \lim_{N\rightarrow\infty}  \frac{w_1(N)}{N}  & \;=\;
    y^*[\xx,1-\yy-2\epsilon,c]   \;=:\; y_1^* \hspace{8mm}\hbox{and}  \\
    \lim_{N\rightarrow\infty}  \frac{w_2(N)}{N}  & \;=\;
    y^*[1-\xx-2\epsilon,\yy,c]   \;=:\; y_2^* \,. 
\end{align*}
Lemma \ref{lem-fmax} assures us that $y^*_1<\xx\wedge(1-\yy-2\epsilon)$ and
$y^*_2<(1-\xx-2\epsilon)\wedge \delta$, and therefore
Equation (\ref{eq.win}) holds for all sufficiently large $N$ (where $I$ is interpreted 
to be $I_N$, etc.).
Using these sequences in Proposition \ref{prop-constr} and invoking
Lemma \ref{lem-Dec} and Equations (\ref{eq.fmax}) and (\ref{eq.tauhatcond}), 
we see that the 
$N^{th}$ root of the right hand side of Equation (\ref{eq.FDDS}) converges to 
\begin{align}
    \nonumber
     & 2^{\xx+1-\yy-2\epsilon}g(\xx,1-\yy-2\epsilon;c)\,
    g(1-\yy-2\epsilon,\xx;c)  \,c^{y_1^*}      \\
   \label{eq.gggglb}
     & \hspace{3mm}\times 
      2^{1-\xx+\yy-2\epsilon}g(1-\xx-2\epsilon,\yy;c)\,
    g(\yy,1-\xx-2\epsilon;c)  \,c^{y_2^*}     
     \times 
      c^{1-y_1^*-y_2^*}  \,.
\end{align}
Thus
Equation (\ref{eq.gggglb})  is a lower bound for 
$\liminf_{N\rightarrow\infty}|{\cal F}^*_{N}(I_N,J_N;\tau)|^{1/N}$
for all sufficiently small positive $\epsilon$.
Now let $\epsilon$ decrease to 0. By the continuity of $g$, the
expression of Equation (\ref{eq.gggglb}) converges
to $G(\xx,1-\yy;c)$.  This proves the proposition.
\hfill $\Box$

\medskip
\noindent
\textbf{Proof of Proposition \ref{prop-constr}:}
Fix $N$, $I$, $J$, $A$, $w_1$ and $w_2$ as specified.
We shall prove the proposition by constructing an injection from 
${\cal D}$ into ${\cal F}^*_{N}(I,J;\tau)$, where 
\begin{align*}
   {\cal D} \;& =  \;
        \Dec^{*A}(w_1;I-1,N-2A-J)  \times
        \Dec^{*A}(w_2;N-2A-I,J-1)    \\
    & \hspace{28mm}
       \times \, {\cal S}^{*A}_{N-w_1-w_2-1}(\hat{\tau}). 
\end{align*}
Consider $({\cal B}_1,{\cal B}_2,\phi)\in {\cal D}$ (that is,
${\cal B}_1$ is one of the  $w_1$-element sets in $\Dec^{*A}(w_1;I-1,N-2A-J)$, 
and so on).
Let $\Psi\equiv \Psi({\cal B}_1,{\cal B}_2)$ be the $(w_1+w_2+1)$-element 
decreasing set defined by
\[   
    \left({\cal B}_1+(0,J)\right)  \,\cup \, \{(I,J)\}  \,\cup \,
    \left({\cal B}_2+(I,0)\right)
\]
(where ${\cal B}+(x,y)$ denotes translation of the set ${\cal B}$ by
the vector $(x,y)$).  Thus $\Psi$ is a 
decreasing subset of $[1,N-2A]\times [1,N-2A]$ that contains $(I,J)$.  

We claim that 
\begin{equation}
   \label{eq.claim1}
  y\;<\;N-x-A \hspace{5mm}\hbox{for every $(x,y)\in \Psi$.}
\end{equation}
For $(x,y)=(I,J)$, this follows from our assumption $J<N-I-2A$.  
For $(x,y)$ in ${\cal B}_1+(0,J)$, we have $(x,y-J)\in {\cal B}_1$
and hence 
\[  
  y-J \;<\; (N-2A-J)\,-\,x\,\frac{N-2A-J}{I-1}\,+\,A
      \;<\, N-A-J-x
\]
(using $I<N-2A-J$), which verifies the claim in this case.  
A similar argument works if $(x,y)\in {\cal B}_2+(I,0)$.
Therefore the claim (\ref{eq.claim1})  is true.

Given $\Psi$ and a permutation 
$\phi\in {\cal S}^{*A}_{N-w_1-w_2-1}(\hat{\tau})$, we shall define a
permutation $\sigma\in {\cal S}_N$ such that $\Psi$ is contained 
in the graph of $\sigma$ (i.e., $y=\sigma_x$ whenever $(x,y)\in \Psi$)
and $\patt(\sigma\setminus\Psi)=\phi$.
Let $w=w_1+w_2+1$, and write the elements of $\Psi$ as 
$(x(\ell),y(\ell))$ ($\ell=1,\ldots,w$) with $x(\ell)$ increasing in $\ell$
and $y(\ell)$ decreasing in $\ell$.  Define the functions
$\Gamma_x$ and $\Gamma_y$ from $\{1,\ldots,N-w\}$ into $\{1,\ldots,N\}$ 
as follows.
Writing $x(0)=0$ and $x(w+1)=N+1$, and observing that 
$x(\ell)-\ell$ is decreasing in $\ell$, we define
\begin{align*}
    \Gamma_x(i)  \;=\;  i+m   \hspace{5mm}
     &  \mbox{where $m$ satisfies  $x(m)-m\,<\,i\,\leq \,x(m+1)-(m+1)$}; \\
     &  \mbox{i.e., where $m$ satisfies  $x(m)\,<\,i+m\,< \,x(m+1)$}.
\end{align*}
The possible values for $m$ are $0,1,\ldots,w$.
Analogously, writing $y(0)=N+1$ and $y(w+1)=0$, we define
\begin{align*}
    \Gamma_y(i)  \;=\;  i+n   \hspace{5mm}
     &  \mbox{where $n$ satisfies  $y(w{-}n{+}1)-n+1\,\leq \,i\,< 
    \,y(w{-}n)-n$}; \\
     &  \mbox{i.e., where $n$ satisfies  $y(w{-}n{+}1)\,<\,i+n\,< \,y(w{-}n)$}.
\end{align*}
Again, the possible values for $n$ range from 0 to $w$.
Observe that $\Gamma_x$ (respectively, $\Gamma_y$) is the unique strictly increasing
function from $\{1,\ldots,N-w\}$ to 
$\{1,\ldots,N\}\setminus \{x(1),\ldots,x(w)\}$   (respectively, 
$\{1,\ldots,N\}\setminus \{y(1),\ldots,y(w)\}$). 
Now define $\sigma_1,\ldots,\sigma_N$ by 
\begin{align*}
    \sigma_{x(\ell)} \;& =\; y(\ell)  \hspace{5mm}
     \hbox{for $\ell=1,\ldots,w$},   \\
    \sigma_{\Gamma_x(i)} \;& =\; \Gamma_y(\phi_i)   
                \hspace{5mm} \hbox{for $i=1,\ldots,N-w$}.
\end{align*}
Then it is not hard to see that the string 
$\sigma :=\sigma_1\sigma_2\cdots \sigma_N$
is well defined, that $\sigma$ is a permutation in ${\cal S}_N$
whose graph contains $\Psi$, and that
$\patt(\sigma\setminus\Psi)=\phi$.
See Figure \ref{fig:exampleA}.

\begin{center}
\begin{figure}[htp]
\centering
  \begin{tikzpicture}[scale=0.23]
    \draw (0,0) rectangle (40,40);
    \draw [red] (0,34) rectangle (16,11);
    \draw [red] (34,0) rectangle (18,9);
    \draw [dashed,blue]  (0,40) -- (40,0);
    \draw [very thin, blue] (0,37) -- (37,0);
    \draw [very thin, red]  (16,14) -- (2.1,34);
    \draw [very thin, red]   (20.5,9)  -- (34,3);
    \draw [fill]     (3,28) circle (0.3);
    \draw [fill]     (9,21) circle (0.3);
    \draw [fill]     (12,15) circle (0.3);
    \draw [fill]     (17,10) circle (0.3);   %
    \draw     (17,10) circle (0.8);   
    \draw [fill]     (23,7) circle (0.3);
    \draw [fill]     (30,3) circle (0.3);            
   \draw [fill]  (0,40) circle (0.3);
   \draw [fill]  (1,39) circle (0.3);
   \draw [fill]  (2,38) circle (0.3);
   \draw [fill]  (4,37) circle (0.3);
   \draw [fill]  (5,36) circle (0.3);
   \draw [fill]  (6,35) circle (0.3);
   \draw [fill]  (7,34) circle (0.3);
   \draw [fill]  (8,33) circle (0.3);
   \draw [fill]  (10,32) circle (0.3);
   \draw [fill]  (11,31) circle (0.3);
   \draw [fill]  (13,30) circle (0.3);
   \draw [fill]  (14,29) circle (0.3);
   \draw [fill]  (15,27) circle (0.3);
   \draw [fill]  (16,26) circle (0.3);
   \draw [fill]  (18,25) circle (0.3);
   \draw [fill]  (19,24) circle (0.3);
   \draw [fill]  (20,23) circle (0.3);
   \draw [fill]  (21,22) circle (0.3);
   \draw [fill]  (22,20) circle (0.3);
   \draw [fill]  (24,19) circle (0.3);
   \draw [fill]  (25,18) circle (0.3);
   \draw [fill]  (26,17) circle (0.3);
   \draw [fill]  (27,16) circle (0.3);
   \draw [fill]  (28,14) circle (0.3);
   \draw [fill]  (29,13) circle (0.3);
   \draw [fill]  (31,12) circle (0.3);
   \draw [fill]  (32,11) circle (0.3);
   \draw [fill]  (33,9) circle (0.3);
   \draw [fill]  (34,8) circle (0.3);
   \draw [fill]  (35,6) circle (0.3);
   \draw [fill]  (36,5) circle (0.3);
   \draw [fill]  (37,4) circle (0.3);
   \draw [fill]  (38,2) circle (0.3);
   \draw [fill]  (39,1) circle (0.3);
   \draw [fill]  (40,0) circle (0.3);
    \draw [dashed] (3,40) -- (3,28) -- (40,28);
     \draw [dashed] (9,40) -- (9,21) -- (40,21);
    \draw [dashed] (12,40) -- (12,15) -- (40,15);
     \draw [dashed] (17,40) -- (17,10) -- (40,10);  
    \draw [dashed] (23,40) -- (23,7) -- (40,7);
     \draw [dashed] (30,40) -- (30,3) -- (40,3);
      \node [below] at (0,0) {$1$};
            \node [below] at (40,0) {$N$};
      \node [below] at (34,0) {$N{-}2A$};
       	\node [left] at (0,40) {$N$};
	\node [left] at (0,34) {$N{-}2A$};
    \node [above] at (3,40) {$x(1)$};	
    \node [above] at (30,40) {$x(w)$};
    \node [above] at (17,40) {$I$};
    \node [right] at (40,28)  {$y(1)$};
    \node [right] at (40,3) {$y(w)$};
    \node [right] at (40,10) {$J$};
  \end{tikzpicture}
  \caption{An example of the permutation $\sigma$ constructed in the proof of 
  Proposition \ref{prop-constr}, in which $N=41$, $w_1=3$, $w_2=2$, $w=6$, and $A=3$, and
  the permutation $\phi$ is the decreasing permutation of length $N-w$.
  The circled black dot is at $(I,J)$.  The dashed blue line is the diagonal of $[1,N]^2$.
  The two red rectangles enclose ${\cal B}_1+(0,J)$ and ${\cal B}_2+(I,0)$.  The sloped red 
  line segment within each red rectangle is drawn $A$ units above the diagonal of the rectangle.
  No point of $\Psi$ is above a sloped red line segment.
  The solid blue line is the line $y=N-x-A$, which partitions the graph of $\sigma$ as 
  described in the Key Claim in the proof.  
  The two sloped red line segments lie below the solid blue line.
  Observe that $I=x(w_1+1)$ and $J=y(w_1+1)$. }
	\label{fig:exampleA}
  \end{figure}
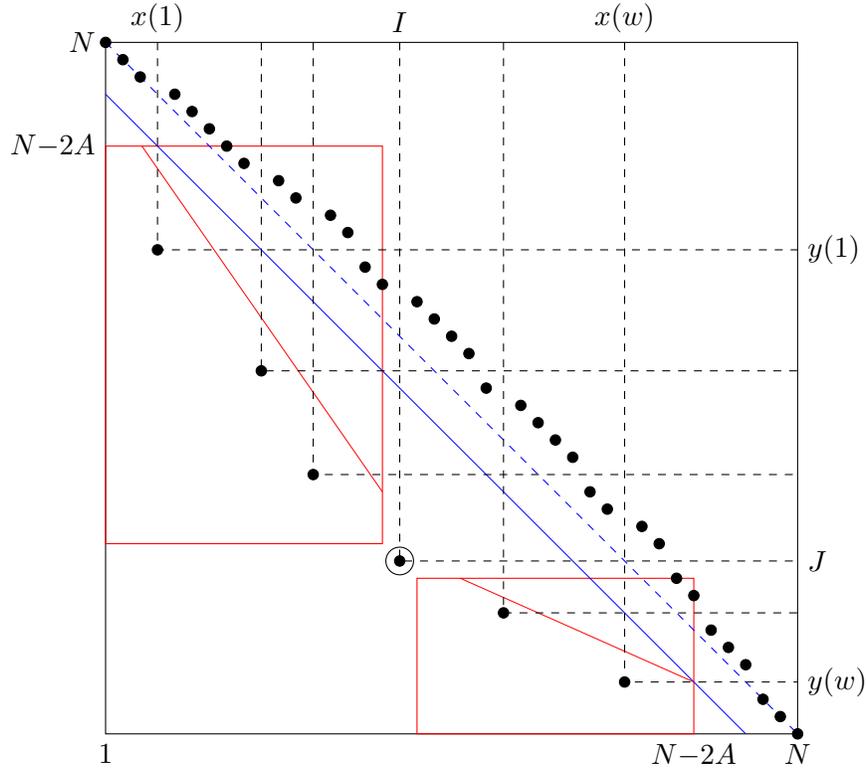
  \end{center}

The proof of the proposition is based on the following claim.
Let $\Psi_x=\{x(1),x(2),\ldots,x(w)\}$.   

\smallskip
\textbf{Key Claim:}  \textit{We have $\sigma_j<N-j-A$ for every $j\in \Psi_x$,
and $\sigma_j>N-j-A$ for every $j\not\in \Psi_x$.}

\smallskip
Once the Key Claim is proven, we proceed as follows.  
The Key Claim implies that $\Psi\subset{\cal M}(\sigma)$  (recall Equation (\ref{eq.defL})).
Therefore, since $Patt(\sigma\setminus\Psi)$ avoids $\hat{\tau}$, so does 
$Patt(\sigma\setminus {\cal M}(\sigma))$.   
Hence, by Lemma \ref{lem-minusL}(\textit{ii}), $\sigma$ avoids $\tau$. 
It follows that $\sigma\in {\cal F}^*_{N}(I,J;\tau)$.
Consequently, writing $Q({\cal B}_1,{\cal B}_2,\phi)\,=\,\sigma$,
we have defined a function $Q:{\cal D}\rightarrow {\cal F}^*_{N}(I,J;\tau)$.
To see that the function $Q$ is one-to-one, suppose
$Q({\cal B}_1,{\cal B}_2,\phi)=\sigma$.
Since $\Psi$ is contained in the graph of $\sigma$, the Key Claim
shows that $\Psi({\cal B}_1,{\cal B}_2)$   
is uniquely determined by $\sigma$, as is $\phi$.
Finally, since $(I,J)$ is specified,
${\cal B}_1$ and ${\cal B}_2$ are determined by $\Psi({\cal B}_1,{\cal B}_2)$.
Hence $Q$ is one-to-one, and the proposition follows.

It only remains to prove the Key Claim.  For $j\in \Psi_x$, say $j=x(\ell)$,
we have $\sigma_j = y(\ell)$, and the assertion of the Key Claim follows from 
Equation (\ref{eq.claim1}).
Now suppose $j\not\in\Psi_x$.  Then for some $i\in [1,N-w]$ we have
$j=\Gamma_x(i)$ and $\sigma_j=\Gamma_y(\phi_i)$.  
Since $\phi\in {\cal S}_{N-w}^{*A}(\hat{\tau})$, we know that $\phi_i>(N-w)-i-A$.
Following the notation in the definitions of $\Gamma_x$ and $\Gamma_y$,
let $m=\Gamma_x(i)-i$ and $n=\Gamma_y(\phi_i)-\phi_i$.
Then $x(m)\,<\,i+m\,<\,x(m+1)$ and 
$y(w{-}n{+}1)\,<\,\phi_i+n\,<\,y(w{-}n)$.
Also, we have
\begin{align*}
   \sigma_j \; & = \;  \phi_i+n   \\
           & > \;  (N-w)-i-A +n   \\
           & = \;  N-w-(j-m) -A +n   \,.
\end{align*}
Thus, to show $\sigma_j>N-j-A$, as required for proving the Key Claim, we need to 
show that $m\geq w-n$.

Assume that $m\geq w-n$ is false, i.e.\ that 
$m+1\leq w-n$.  Since $y(\ell)\geq y(\ell+1)+1$ for every $\ell$, we
see that
\begin{equation}
   \nonumber
      y(m+1) \;\geq \; y(w-n) \,+\,(w-n)-(m+1)  \,.
\end{equation}
Using this inequality and those of the preceding paragraph, we obtain
\begin{align*}
    N-w-A  \;& < \;  \phi_i+i    \\
     & \leq \;   y(w-n)-n-1 +x(m+1)-m-1   \\
     & \leq \;  [y(m+1) -w+n+m+1]  -n+x(m+1)-m-2
   \\
    & \leq \;  N-A-w -1
        \hspace{25mm}\hbox{(by (\ref{eq.claim1}))}   
\end{align*}
which is a contradiction.  Therefore $m\geq w-n$.
This proves  the Key Claim, and hence the proposition.
\hfill $\Box$
       
\medskip
\noindent
\textbf{Proof of Lemma \ref{lem-Dec}:}
For positive integers $w$ and $M$, let $\Seq(w;M)$ be the set of
all $w$-element subsets of $\{1,2,\ldots,M\}$.  We shall write 
a member of $\Seq(w;M)$ as a $w$-element vector with
the entries in increasing order:  $\vec{x}=(x(1),x(2),\cdots,x(w))$, with
$x(1)<\cdots < x(w)$.
Then there is a natural bijection 
$\Theta:\Seq(w;M_1)\times \Seq(w;M_2)\rightarrow\Dec(w;M_1,M_2)$ via
\[
   \Theta(\vec{x},\vec{z}) \;=\;  \{(x(1),z(w)),(x(2),z(w-1)),\cdots,
     (x(w),z(1))\} \,.
\]
In particular, we have
\begin{equation}
   \label{eq.seqcard}
    |\Dec(w;M_1,M_2)|  \;=\; |\Seq(w;M_1)| \,|\Seq(w;M_2)|  \;=\;
       \binom{M_1}{w} \,\binom{M_2}{w}  \,.
\end{equation}
Applying Lemma \ref{lem-Stir} to Equation (\ref{eq.seqcard}) proves
Equation  (\ref{eq.lemDec2}).  Equation (\ref{eq.lemDec2A}) will 
follow immediately once we have proven Equation (\ref{eq.lemDec1}).

For positive integers $A$, we now define 
\[    \Seq^{*A}(w;M)  \;=\;  \left\{\vec{x}\in \Seq(w;M):
     \left|x(\ell)-\ell\,\frac{M}{w+1}\right|<A \hbox{ for }\ell=1,\ldots,w
     \right\} \,.
\]
Roughly speaking, 
a $w$-element subset of $\{1,\ldots,M\}$ is in $\Seq^{*A}(w;M)$
if its elements are within distance $A$ of a uniform spacing configuration 
over the interval.  We shall now show the following.
\begin{align*}
   & \textbf{Property I}:  \;\;
  \textrm{If $\vec{x}\in \Seq^{*A}(w;M_1)$ and $\vec{z}\in\Seq^{*A}(w;M_2)$, 
  then}   \\
   & \textrm{ $\Theta(\vec{x},\vec{z})\in \Dec^{*B}(w;M_1,M_2)$, where
  $B=A\left(1+\frac{M_2}{M_1}\right)$.}
\end{align*}
Property I says that if $\vec{x}$ and $\vec{z}$ are close 
to being uniformly spaced on their intervals, then $\Theta(\vec{x},\vec{z})$
is close to the diagonal of its rectangle.
To prove Property I, consider $\vec{x}$ and $\vec{z}$ as specified.
Then a generic point of $\Theta(\vec{x},\vec{z})$, $(x(\ell),z(w+1-\ell))$,
satisfies
\begin{align*}
   &\left|z(w+1-\ell)-\left(M_2  -x(\ell)\,\frac{M_2}{M_1}\right)\right|   \\
     & \hspace{8mm} \leq  \;  \left|z(w+1-\ell)-(w+1-\ell)\,\frac{M_2}{w+1}\right| 
     \,+\,  \frac{M_2}{M_1}\left|x(\ell)-\ell\,\frac{M_1}{w+1}\right|   \\
     & \hspace{8mm} < \;A +\frac{M_2}{M_1}\,A \,.
\end{align*}
This proves Property I.
Now, Property I implies that $|\Dec^{*B}(w;M_1,M_2)|\,\geq \,
|\Seq^{*A}(w;M_1)|\,|\Seq^{*A}(w;M_2)|$.  Recalling Equation (\ref{eq.seqcard}), we see
that Equation (\ref{eq.lemDec1}) will follow if we can prove
\begin{align*}
   \textbf{Property II}:  \;\;& \lim_{N\rightarrow\infty} 
   \frac{|\Seq^{*A_N}(w(N);N)|}{\binom{N}{w(N)} }  \,=\,1   
     \hspace{5mm} \hbox{ whenever }   \\
  & \lim_{N\rightarrow\infty}\frac{w(N)}{N}\,=:\,\theta\,\in\,(0,1)
     \hspace{3mm}\hbox{ and }  \hspace{3mm}
   \lim_{N\rightarrow\infty}\frac{A_N}{N}\,=:\,\epsilon\,>\,0.
\end{align*}

We shall prove Property II by converting it into a probabilistic statement.
Let $p\in (0,1)$.  Let $G_1,G_2,\ldots$ be a sequence of independent
random variables having the  geometric distribution with parameter $p$;
that is, $\Pr(G_i=\ell)=p(1-p)^{\ell-1}$ for $\ell=1,2,\ldots$.
Next, let $T_i=G_1+G_2+\cdots+G_i$ for each $i$.
These random variables have negative binomial distributions
\begin{equation}
   \label{eq.negbin}
    \Pr(T_{j+1}=\ell+1)  \;=\;  \binom{\ell}{j}p^{j+1}(1-p)^{\ell-j}  
    \hspace{5mm}\hbox{for $\ell \geq j$}.
\end{equation}
Moreover, for any $\vec{x}\in \Seq(w;N)$ (writing $x(0)=0$ and $x(w+1)=N+1$),
\begin{align}
   \nonumber
   \Pr(T_{\ell}=x(\ell) \hbox{ for }\ell=1,\ldots,w \,|\, T_{w+1}=N{+}1)
    \; & =\;  \frac{ \prod_{\ell=1}^{w+1} p(1-p)^{x(\ell)-x(\ell-1)-1} }{
         \binom{N}{w}p^{w+1}(1-p)^{N-w} }
   \\  \label{eq.geomjoint}
     & = \;   \binom{N}{w}^{-1}  \,.
\end{align}  
Equation (\ref{eq.geomjoint}) says that
\textit{the conditional distribution of $(T_1,\ldots,T_w)$ given 
that $T_{w+1}=N+1$ is precisely the uniform distribution on $\Seq(w;N)$.}
This assertion is true for any $p$.  Let us now fix $p=(w+1)/N$; we shall
soon see why this is a convenient choice.

By Equation (\ref{eq.geomjoint}),  
\begin{equation}
    \nonumber
   \frac{ |\Seq^{*A}(w;N)|}{\binom{N}{w}}  \; = \;
     \Pr\left( |T_{\ell}-\ell/p| < A \hbox{ for }l=1,\ldots,w \,|\, 
       T_{w+1}=N+1  \right)   \,.
\end{equation}
and therefore
\begin{equation}
    \label{eq.ratioprob}
   0\;\leq \; 1 \,-\, \frac{ |\Seq^{*A}(w;N)|}{\binom{N}{w}}  
   \; \leq  \; \frac{ \Pr(\max_{\ell=1,\ldots,w}|T_{\ell}-\ell/p| \geq  A )}{
        \Pr(T_{w+1}=N+1)} .
\end{equation} 
It is straightforward to derive the asymptotic behaviour 
$\Pr(T_{w+1}=N+1)$ using Stirling's Formula 
$m!\sim \sqrt{2\pi m}(m/e)^m$ and $p=(w+1)/N$, with 
$w=w(N)\sim \theta N$, as follows.
\begin{align}
    \nonumber
   \Pr(T_{w+1}=N+1)  \; &=\;
     \frac{N!}{w!(N-w)!}  \frac{(w+1)^{w+1}(N-w-1)^{N-w}}{N^{N+1} }    \\
    \nonumber
    & \sim \;  \frac{\sqrt{2\pi N}}{\sqrt{2\pi w}\sqrt{2\pi(N-w)}}
      \left(\frac{w{+}1}{w}\right)^{w+1} \frac{w}{N}
      \left(\frac{N{-}w{-}1}{N{-}w}\right)^{N-w}   \\
   & \sim \;   \frac{\sqrt{\theta}}{\sqrt{2\pi (1-\theta)N} } \,.
     \label{eq.denom}
\end{align}
For the numerator of the right-hand side of Equation (\ref{eq.ratioprob}),
we use Kolmogorov's Inequality \cite{Feller}, 
along with the property that the random 
variables $G_i$ have mean $1/p$ and variance $(1-p)/p^2$:
\begin{align}
   \nonumber
  \Pr\left(\max_{\ell=1,\ldots,w}|T_{\ell}-\ell/p| \geq  A \right)
    \;& \leq \;  \frac{\textrm{Var}(T_w)}{A^2}    \\
     & \sim  \;   \frac{N(1-\theta)/\theta^2}{N^2\epsilon^2}   
    \label{eq.numer}
\end{align}
Applying  Equations (\ref{eq.denom}) and (\ref{eq.numer}) 
to Equation (\ref{eq.ratioprob}) proves Property II.
This completes the proof of Lemma \ref{lem-Dec}.
\hfill $\Box$

\section{Conclusion}
   \label{sec-conc}

Recalling Remark \ref{rem.LDgen}(b), we see that Theorem \ref{thm.LDgen} 
follows immediately from Propositions \ref{prop.Fnostar} and \ref{prop.Fstar2}.

We now show that Equation (\ref{eq.limF}) of Theorem \ref{thm.LD} follows from 
Theorem \ref{thm.LDgen} by induction.   
Remark \ref{rem.LDgen}(a) tells us that we can apply Theorem \ref{thm.LDgen}
when $\tau$ is $1\oslash \mu_3$, 
which shows that Equation (\ref{eq.limF}) holds for $k=4$.  
Now assume that Equation (\ref{eq.limF}) is 
true for a given $k\geq 4$.  Lemma \ref{lem.Gmono} and Remark \ref{rem-diag}
prove that $G(\xx,1-\yy;  (k-2)^2) \,<\, (k-1)^2$ whenever
$\xx < 1-\yy$.  This means that Equation (\ref{eq.limF})   
implies Equation (\ref{eq.tauhatcond}) when $\hat{\tau}$ is $\mu_k$,  
using 
\[   {\cal S}_N(\mu_k)\setminus {\cal S}_N^{*N\epsilon}(\mu_k) 
   \; \;\subset \;   \bigcup_{i,j\,:\; j\,\leq \,N-i-N\epsilon} {\cal F}_N(i,j;\mu_k)
\]
and a compactness argument as in the proof of Proposition \ref{prop.Fnostar}.  
Hence Equation (\ref{eq.limFgen}) holds when $\tau$ is $\mu_{k+1}$,
in which case $L(\hat{\tau})$ equals $(k-1)^2$.  This says that 
Equation (\ref{eq.limF}) %
holds with $k$ replaced by $k+1$.
This completes the induction, showing that  Equation (\ref{eq.limF}) holds for every $k\geq 4$.

Finally we shall prove Equation (\ref{eq.limF2}) for $k\geq 4$ and $1\leq \ell\leq k-2$.  
The proof of Proposition 2.3 in \cite{backelin} shows that there is a bijection from 
${\cal S}_N(1\ldots \ell(\ell{+}1)\cdots(k{-}1)k)$
to ${\cal S}_N(1\ldots \ell k(k{-}1)\ldots(\ell+1))$ that preserves all the left-to-right minima of
each permutation.  (To see this, observe that  when $A=J_{\ell}$ in the proof of \cite{backelin}, 
each right-to-left minimum and everything below it and to its right are all coloured blue, and
hence are unchanged by the bijection $\alpha$.)
It follows that 
\[     {\cal F}^*_N(I,J;\lambda_{k,\ell})   \;=\;   {\cal F}^*_N(I,J;\mu_k)    \]
always holds.  Using this and our Proposition \ref{prop.Fstar2} with $\tau=\mu_k$, we obtain
\begin{equation}
   \label{eq.limFstar4}
   \liminf_{N\rightarrow\infty}|{\cal F}^*_{N}(I_N,J_N;\lambda_{k,\ell})|^{1/N}  \;\geq \;
     G(\xx,1-\yy; (k-2)^2).
\end{equation}
Next, by Proposition \ref{prop.Fnostar} with $\tau=\lambda_{k,\ell}$, we obtain
\begin{eqnarray}
   \nonumber
   \limsup_{N\rightarrow\infty}|{\cal F}_{N}(I_N,J_N;\lambda_{k,\ell})|^{1/N}  &\leq &
     G(\xx,1-\yy;L(\lambda_{k-1,\ell-1}))    \\
     \label{eq.limFnostar2}
      &= & G(\xx,1-\yy;(k-2)^2)  .
\end{eqnarray}
Equations (\ref{eq.limFstar4}) and (\ref{eq.limFnostar2}) together imply Equation (\ref{eq.limF2}).
This completes the proof of  Theorem \ref{thm.LD}.

\section{Acknowledgments}
The research of N. Madras was supported in part by a
Discovery Grant from the Natural  Sciences and
Engineering Research Council of Canada.
Part of this work was done while N. Madras was visiting the Fields Institute for Research in Mathematical Sciences.    
L. Pehlivan would like to thank to the Department of Mathematics and Statistics at Dalhousie University
for their hospitality while she was working on the paper.
The authors thank Erik Slivken for informative discussions.

\bibliography{2016-05-24pattern}
\end{document}